%% file: agt-5-70.tex
\documentclass{gtart_h}
\input agtout

\lognumber{70}
\volumenumber{5}
\volumeyear{2005}
\papernumber{70}
\pagenumbers{1719}{1732}
\received{17 May 2005} 
\accepted{2 December 2005}
\published{17 December 2005}

\usepackage{amsmath,amssymb,graphicx}

\def\fref#1{\hyperlink{#1anchor}{\ref*{#1}}}
\def\figref#1{\hyperlink{#1anchor}{Figure~\ref*{#1}}}
\def\anchor#1{\noindent\hypertarget{#1anchor}{\smash{$\phantom{99}$}}}%\newline}

\newtheorem{thm}{Theorem}[section] 
\newtheorem{lemma}[thm]{Lemma}
\newtheorem{proposition}[thm]{Proposition}
 
\newtheorem{corollary}[thm]{Corollary}
\newtheorem{conjecture}[thm]{Conjecture}
\newtheorem{cor}{Corollary}[section]

\theoremstyle{definition}

\newtheorem{defi}[cor]{Definition}

\def\Wh{{\rm Wh}}
\def\Bor{{\rm Bor}}
\def\Free{{\rm Free}}

\newcommand{\cover}{
$$\begin{picture}(140,150)  \small
    \put(-115,13)       {\includegraphics{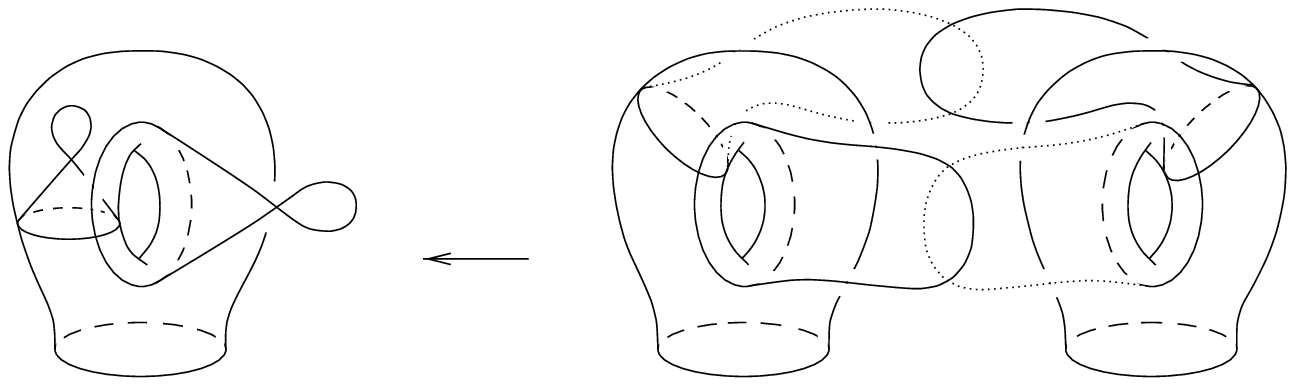}}
   \put(-47,20)       {$\gamma$}
   \put(-120,31)      {$S^c$}
   \put(-16,73)       {$A$}
   \put(-93,94)       {$B$}
   \put(55,20)        {$\overline{\gamma}_1$}
   \put(138,26)       {$\overline{A}_1$}
   \put(125,125)      {$\overline{B}_1$}
   \put(175,125)      {$\overline{B}_2$}
   \put(162,26)       {$\overline{A}_2$}
   \put(247,20)       {$\overline{\gamma}_2$}
\end{picture}$$}

\newcommand{\double}{
$$\begin{picture}(160,195) \small
    \put(-105,5)        {\includegraphics{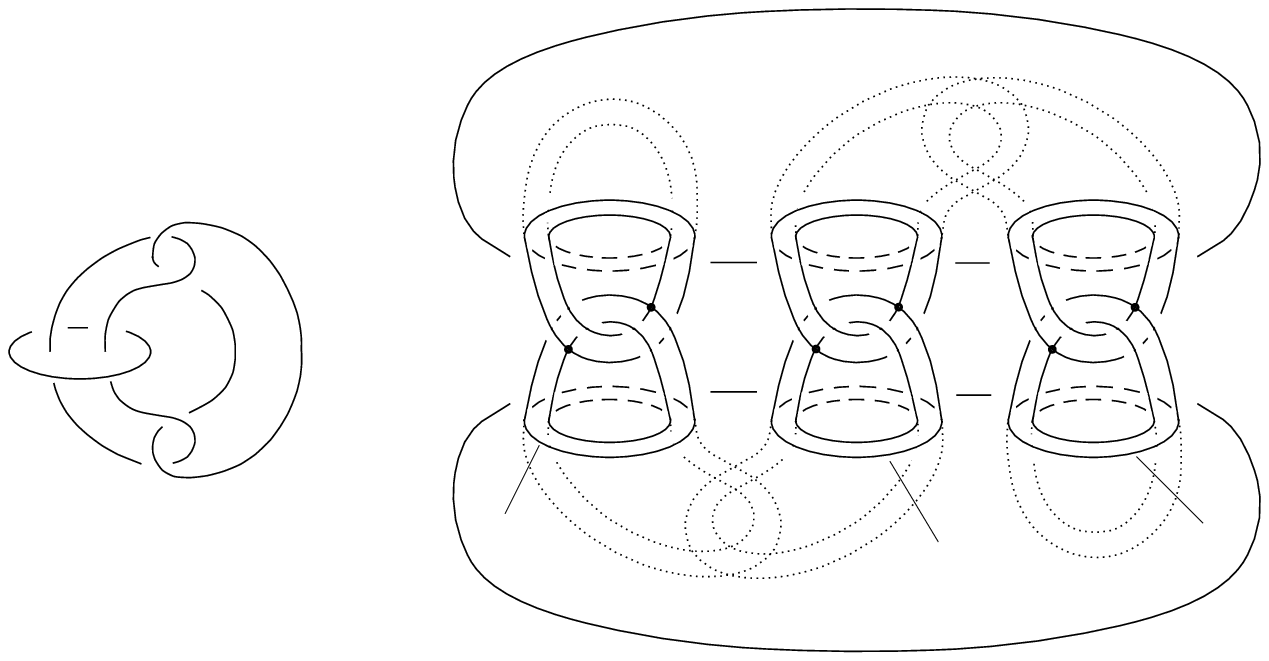}}
    \put(-95,40)     {$\Bor=(l_1,l_2,l_3)$}
    \put(7,60)       {$B^4$}
    \put(7,120)      {$\overline B^4$}
    \put(33,35)      {$l_1,l'_1$}
    \put(157,26)     {$l_2,l'_2$}
    \put(230,32)     {$l_3,l'_3$}
\end{picture}$$}

\newcommand{\plumbing}{
$$\begin{picture}(140,90) 
    \put(-50,5)        {\includegraphics{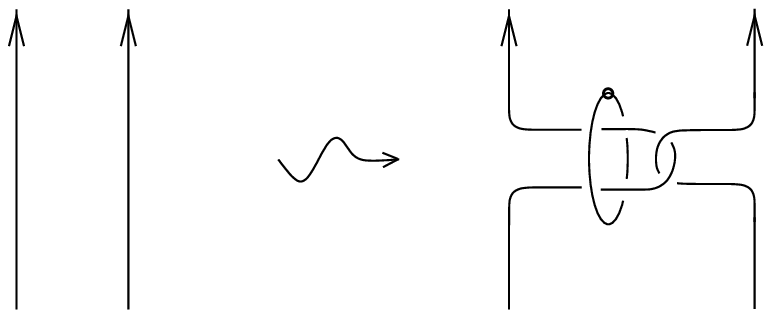}}
\end{picture}$$}

\newcommand{\graphcover}{
$$\begin{picture}(140,90) \small
    \put(-90,5)        {\includegraphics{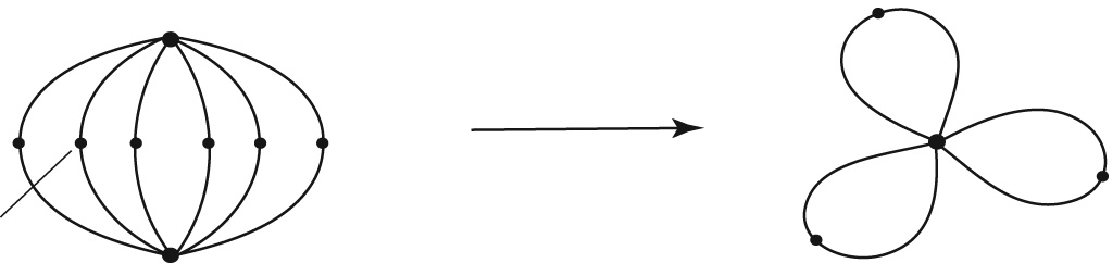}}
    \put(-101,39)  {$q_1$}
    \put(-100,13)   {$q'_1$}
    \put(153,83)     {$p_1$}
    \put(133,5)      {$p_2$}
    \put(234,28)     {$p_3$}
    \put(9,39)    {$q'_3$}
    \put(-95,65)  {$\Gamma$}
\end{picture}$$}

\newcommand{\equivalence}{
$$\begin{picture}(140,90) \small
    \put(-115,5)        {\includegraphics{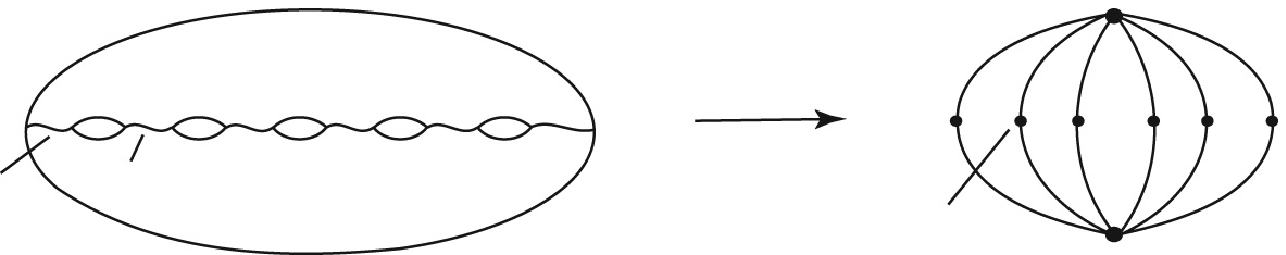}}
    \put(-38,15)    {$N_0$}
    \put(-38,58)    {$N_1$}
    \put(-118,19)   {$Y_1$}
    \put(-87,21)    {$Y'_1$}
    \put(25,0)       {$N$}
    \put(240,5)     {$\Gamma$}
    \put(102,49)    {$f$}
    \put(147,42)    {$q_1$}
    \put(149,11)    {$q'_1$}
\end{picture}$$}

\newcommand{\mathcalc}{
$$\begin{picture}(140,190) \small
    \put(-107,5)        {\includegraphics[width=.95\hsize]{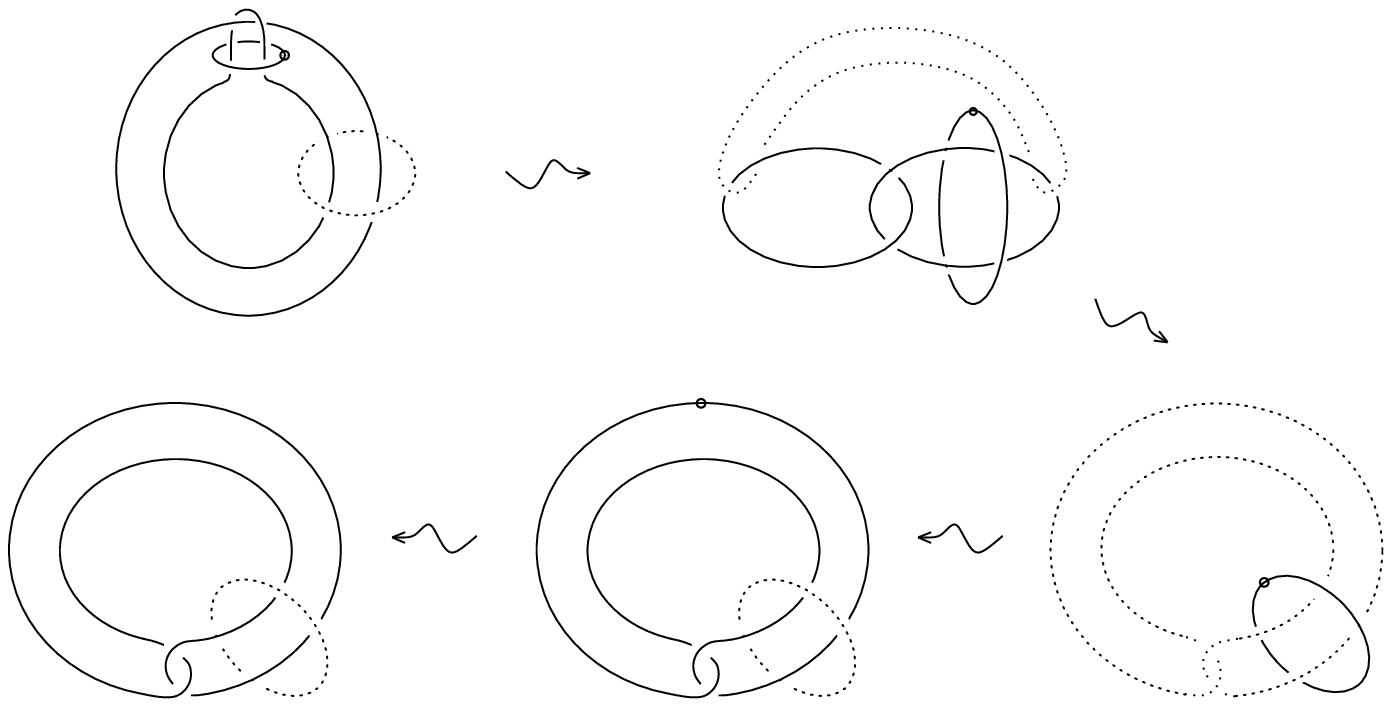}}
    \put(-45,87)      {$0$}
\end{picture}$$}

\begin{document}

\title{Surgery and involutions on 4-manifolds}
\author{Vyacheslav S. Krushkal}

\address{Department of Mathematics, University of 
Virginia\\Charlottesville, VA 22904, USA}
\email{krushkal@virginia.edu}

\begin{abstract}   % type your abstract below
We prove that the canonical $4$-dimensional surgery problems can be
solved after passing to a double cover. This contrasts the
long-standing conjecture about the validity of the topological
surgery theorem for arbitrary fundamental groups (without passing to
a cover). As a corollary, the surgery conjecture is reformulated in
terms of the existence of free involutions on a certain class of
$4$-manifolds. We consider this question  and analyze its relation
to the $A,B$-slice problem.
\end{abstract}

\asciiabstract{%
We prove that the canonical 4-dimensional surgery problems can be
solved after passing to a double cover. This contrasts the
long-standing conjecture about the validity of the topological surgery
theorem for arbitrary fundamental groups (without passing to a
cover). As a corollary, the surgery conjecture is reformulated in
terms of the existence of free involutions on a certain class of
4-manifolds. We consider this question and analyze its relation to the
A,B-slice problem.}


\maketitle  %%% Makes a short header for AGT and GTM articles

\section{Introduction}

The geometric classification techniques --- surgery and the
s-cobordism theorem --- are known to hold in the topological category
in dimension $4$ for a class of fundamental groups which includes
the groups of subexponential growth \cite{F0}, \cite{FT1},
\cite{KQ}, and are conjectured to fail in general \cite{F}. The
unrestricted surgery theorem is known to be equivalent to the
existence of a certain family of canonical $4$-manifolds with free
fundamental group. Recall the precise conjecture concerning these
canonical surgery problems \cite{F}:

\begin{conjecture} \label{conjecture} 
The untwisted Whitehead double of the Borromean Rings, $\Wh(\Bor)$, is
not a freely topologically slice link.
\end{conjecture}

In this statement the additional ``free'' requirement is that the
complement of the slices in the $4$-ball has free fundamental group
generated by the meridians to the link components. (The slicing
problem is open without this extra condition as well.) Considering
the slice complement, the conjecture is seen to be equivalent to the
statement that there does not exist a topological $4$-manifold $M$,
homotopy equivalent to $\vee^3 S^1$, and whose boundary is
homeomorphic to the zero-framed surgery on the Whitehead double of
the Borromean rings: $\partial M\cong {\mathcal{S}}^0(\Wh(\Bor))$. In
contrast, here we show that there exists a double cover of this
hypothetical manifold $M$:

\begin{thm} \label{double cover} 
There exists a smooth $4$-manifold $N$ homotopy equivalent to a
double cover of $\vee^3 S^1$ with $\partial N$ homeomorphic to the
corresponding double cover of ${\mathcal{S}}^{0}(\Wh(\Bor))$. The
analagous double covers exist for all generalized Borromean rings.
\end{thm}

Here the generalized Borromean rings are a family of links obtained
from the Hopf link by iterated ramified Bing doubling. The Borromean
rings are the simplest representative of this class of links, and
the slicing problem for the Whitehead doubles of such links provides
a set of canonical surgery problems.

It is interesting to note that the double covers of the
(hypothetical) canonical $4$-manifolds, constructed in theorem
\ref{double cover}, are {\em smooth}. The surgery conjecture is
known to fail in the smooth category even in the simply-connected
case \cite{D}, thus there does not exist a smooth free involution on
$N$ --- corresponding to at least one of the generalized Borromean
rings --- extending the obvious involution on the boundary $\partial
N$. In the topological category we have the following new
reformulation of the surgery conjecture.

\begin{corollary}  \label{involution} 
Suppose the topological $4$-dimensional surgery and $5$-dimen-
sional s-cobordism conjectures hold for free groups. Then there
exists a free involution on each manifold $N$ constructed in theorem
\ref{double cover}, extending the given involution on the boundary.
Conversely, suppose such involutions exist. Then the surgery
conjecture holds for all fundamental groups.
\end{corollary}

To have a statement equivalent just to surgery (independent of the
s-cobordism conjecture), one needs to consider involutions on the
class of all $4$-manifolds with the given boundary and homotopy
equivalent to those constructed in theorem \ref{double cover}. To
prove the corollary, note that if a required involution existed then
the quotients provide solutions to the  canonical surgery problems.
Conversely, the validity of the surgery conjecture implies that a
surgery kernel (a direct sum of hyperbolic pairs $\left(
\begin{smallmatrix} 0 & 1 \\ 1 & 0 \end{smallmatrix} \right) $ in
${\pi}_2 M$) can be represented by embedded spheres is a
$4$-manifold s-cobordant to $M$, see \cite[12.3]{FQ}. The proof of
theorem \ref{double cover} (see section \ref{canonical problems})
involves a construction of a specific manifold $M_{\Bor}$ surgery on
which would give a canonical manifold as above. Therefore if the
surgery conjecture were true there exists a manifold $M'$
s-cobordant to $M_{\Bor}$ where the surgery kernel is represented by
embedded spheres. Surger them out and consider its double cover
$N'$. If the s-cobordism conjecture holds then $N'$ is homeomorphic
to the manifold $N$ constructed in theorem \ref{double cover}, so
$N$ admits a free involution.\endproof

The proof of theorem \ref{double cover} is given in section
\ref{canonical problems}. Another set of canonical problems is
provided by {\em capped gropes} (they are canonical in the sense
that finding an embedded disk in capped gropes is equivalent to both
the surgery and the s-cobordism theorems.) Section \ref{canonical
problems} also contains a proof analogous to theorem \ref{double
cover} in this context. In section \ref{undoubling}  we start
analyzing the approach to the surgery conjecture provided by the
corollary \ref{involution} above. Its relation to the $A,B$-slice
problem is illustrated by showing that an involution cannot have a
fundamental domain bounded by certain homologically simple
$3$-manifolds. The complexity of this problem is in the interplay
between the topology of $\partial N$ and the homotopy type of $N$ ---
we point out how the analogous question is settled when $N$ is
closed or has a ``simpler'' boundary.

\medskip

{\bf Acknowledgements}\qua I would like to thank Michael
Freedman and Frank Quinn for discussions on the subject. I also
thank the referee for the comments on the exposition of the paper.

This research was partially supported by the NSF and by the Institute for
Advanced Study.

\section{Double cover of the canonical problems}
\label{canonical problems}

In this section we prove Theorem \ref{double cover}. The proof
consists of an explicit construction of a $4$-manifold with the
prescribed boundary, and an observation that the surgery kernel in a
double cover is represented by embedded spheres. Surgering them out
gives a $4$-manifold with the required homotopy type. We also
present an argument for a different set of canonical disk embedding
problems: capped gropes.

{\bf Proof of Theorem \ref{double cover}}\qua Given a link
$L=(l_1,\ldots,l_n)$ in $S^3$, let $L'=(l'_1,\ldots,l'_n)$ denote
its parallel copy, where the components are pushed off with trivial
linking numbers. Consider the $4$-manifold $M_L$ obtained by
attaching zero framed $2$-handles to $B^4$ along the
$(2n)$-component link $L\cup L'$, and introducing a plumbing between
the handles attached to $l_i$ and $l'_i$, for each $i$. The
fundamental group of $M_L$ is the free group $F_{n}$, freely
generated by $n$ loops, each passing exactly once through a plumbing
point. Suppose all linking numbers of $L$ vanish, then the
intersection form on ${\pi}_2(M_L)$ is a direct sum of hyperbolic
planes $\left( \begin{smallmatrix} 0 & 1 \\ 1 & 0 \end{smallmatrix}
\right)$.

The boundary of $M_{L}$ is diffeomorphic to the zero-framed surgery
on $S^3$ along the untwisted Whitehead double of $L$, cf \cite{F1},
\cite{FT1}. For convenience of the reader and since the manifold
$M_L$ (where $L=$ the Borromean rings) plays an important role in
the proof of theorem \ref{double cover}, we provide a proof based on
Kirby calculus \cite{FT1}.

\begin{lemma}\label{calc lemma} 
$\partial M_L\cong {\mathcal{S}}^0(\Wh(L))$, with the isomorphism
carrying the meridians of the $1$-handles to the meridians to
$\Wh(L)$.
\end{lemma}

\begin{figure}[ht]
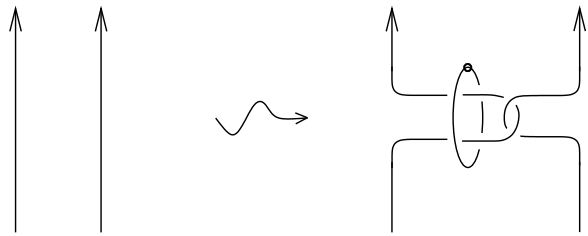
\anchor{plumbing}
\plumbing \caption{Kirby diagram of a (positive) plumbing }
\label{plumbing}
\end{figure}

There is a $\pm$ ambiguity for the clasp of each component (the sign
of the Whitehead doubling corresponds to the sign of the plumbing),
and since it is irrevelevant for our discussion, we consider any
choice of the sign for each component, so $\Wh(L)$ denotes any of the
$2^n$ resulting links. Recall a well-known fact in Kirby calculus:

\begin{proposition} \label{kirby} 
The effect of introducing a $\pm$ plumbing on the underlying Kirby
handle diagram of a handlebody is to introduce a new $1$-handle and
a $\pm$ clasp of the attaching curves of the $2$-handles being
plumbed over the $1$-handle (\figref{plumbing}).
\end{proposition}

\begin{figure}[ht]\anchor{calc}
\mathcalc \caption{Proof of lemma \ref{calc lemma}} \label{calc}
\end{figure}

\proof[Proof of lemma \ref{calc lemma}] This is a calculation in Kirby
calculus in a solid torus, see \figref{calc}. The solid torus is
the complement of the dotted circle in $S^3$. \figref{calc}
displays the positive clasp, of course the negative case is treated
analogously. The first and third arrows are isotopies of ${\mathbb
R}^3$. The second arrow involves replacing a zero-framed $2$-handle
with a $1$-handle (a diffeomorphism of the boundary) and then
cancelling a $1$- and $2$-handle pair. Similarly, the last arrow is
a diffeomorphism of the boundary of the $4$-manifold. \endproof

The construction of $M_L$ is used in \cite{FT1} to show that the
untwisted Whitehead doubles of a certain subclass of homotopically
trivial links are slice. This is done by representing a hyperbolic
basis of ${\pi}_2(M_L)$ by ${\pi}_1$-null transverse pairs of
spheres, which by \cite[Chapter 6]{FQ} are s-cobordant to embedded
pairs. We show that for the generalized Borromean Rings, one can
find {\em smoothly embedded} transverse pairs in a double cover. (Of
course these links are homotopically {\em essential} and this is the
central open case in the surgery conjecture.)

\begin{figure}[ht]
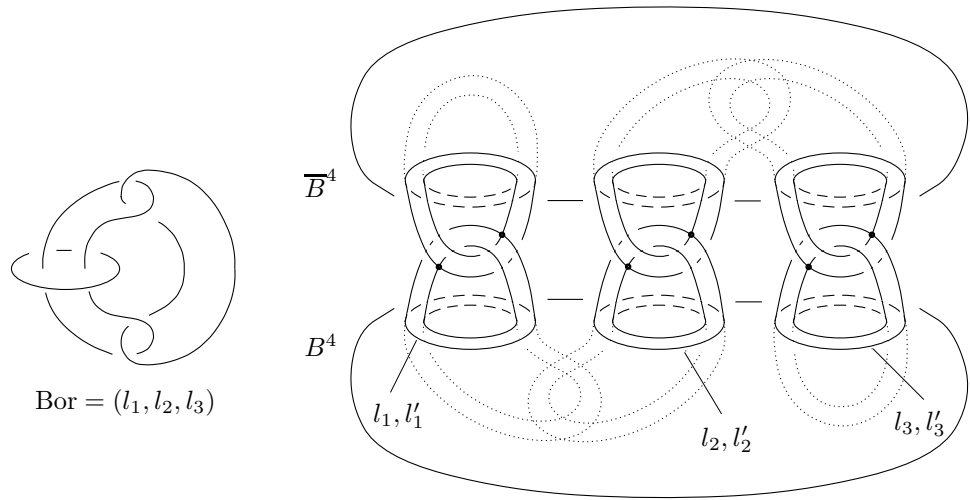
\anchor{double figure}
\double \caption{The Borromean rings, and the double cover
$\overline{M}_{\Bor}$ } \label{double figure}
\end{figure}

\begin{figure}[t]
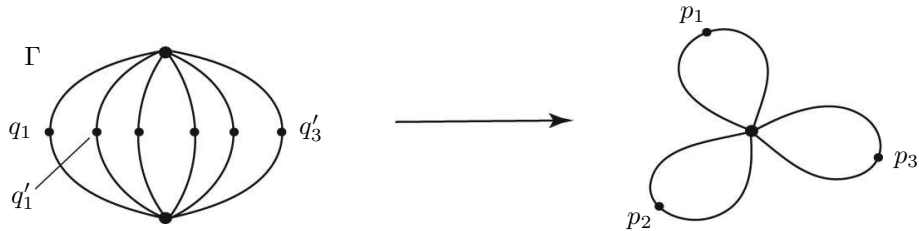
\anchor{graphcover}
\small \graphcover \caption{The double cover
${\Gamma}\longrightarrow S^1\vee S^1\vee S^1$, corresponding to the
homomorphism $\Free_3\longrightarrow {\mathbb Z}/2{\mathbb Z}$,
sending each ``preferred'' generator to the non-trivial element.}
\label{graphcover}
\end{figure}

Consider the Borromean Rings $\Bor=(l_1,l_2,l_3)$, and the
corresponding manifold $M_{\Bor}$, introduced above. Let
$p\co\widetilde{M}_{\Bor}\longrightarrow M_{\Bor}$ be the double cover
induced by the homomorphism
${\pi}_1(M_{\Bor})\longrightarrow{\mathbb{Z}}/2$, which sends each
preferred free generator to the non-trivial element. More
specifically, the $4-$ball and the $2-$handles of $M_{\Bor}$ each
have two lifts in the double cover. Denote the two $4-$balls by
$B^4$, $\overline B^4$, and consider the lifts $\Bor$,
$\overline{\Bor}$ of the Borromean rings and their parallel copies:
$l_i, l'_i$ in $\partial B^4$; $\overline l_i, \overline l'_i$ in
$\overline B^4$, $i=1,2,3$, see \figref{double figure}. The
cover $\widetilde{M}_{\Bor}$ is obtained from
$$(B^4\cup_{\Bor, \Bor'} 2-{\rm handles})\amalg(\overline
B^4\cup_{\overline{\Bor}, \overline{\Bor'}} 2-{\rm handles})$$ by
introducing a plumbing of the handles attached to $l_i, \overline
l'_i$ and also of the handles attached to $l'_i, \overline l_i$, for
each $i$. There is an obvious involution $\tau$ on
$\widetilde{M}_{\Bor}$, with $\widetilde{M}_{\Bor}/{\tau}\cong
M_{\Bor}$, and comparing with \figref{graphcover}, one observes that
this is the double cover corresponding to the required homomorphism
${\pi}_1(M_{\Bor})\longrightarrow{\mathbb{Z}}/2$.

Observe that ${\pi}_2(\widetilde{M}_{\Bor})$ consists of six
hyperbolic planes, and we represent their bases by $2$-spheres as
shown in \figref{double figure}. The special intersection point
in each hyperbolic pair is a lift in the cover of the plumbing point
of the cores of the $2$-handles attached to $B^4$. These cores are
capped off with disks (drawn dotted in the figure) bounded by the
Borromean rings in the two lifts of the $4$-ball, and the extra
intersection points between the spheres are the intersections
between the disks. The key point of the argument is the choice of
the disks bounded by the Borromean rings. Two different ways of
unlinking them are used in the two lifts of the $4-$ball: in one of
them, $B^4$, the components $l_2, l'_2$ intersect $l_1, l'_1$. In
the other one, $\overline B^4$, $\overline l_2, \overline l'_2$
intersect $\overline l_3, \overline l'_3$. To be specific, we
introduce a notation for the spheres formed by the cores of the
$2-$handles capped off with the disks: $S_i, S'_i$, $i=1,2,3$ are
the spheres in $B^4\cup 2-$handles attached to $\Bor\cup \Bor'$. Here
the index reflects the component of the link giving rise to the
sphere. The spheres in the other lift, $\overline B^4\cup 2-$handles
attached along $\overline{\Bor}\cup\overline{\Bor'}$, are denoted by
$\overline S_i$, $\overline S'_i$, $i=1,2,3$. The six hyperbolic
planes in ${\pi}_2(\widetilde M_{\Bor})$ are formed by the pairs of
$2-$spheres $(S_i, \overline S'_i)$ and $(\overline S_i, S'_i)$,
$i=1,2,3$. The spheres in each pair intersect in precisely one
point: the plumbing point of the corresponding $2-$handles.

Note that none of the spheres have self-intersections (this is true
even without taking a cover, in $M_{\Bor}$). Moreover, $\overline
S_1, \overline S'_1, S_3, S'_3$ are embedded disjointly from other
spheres, except for the special intersection (plumbing) points, and
they provide embedded transverse spheres (cf \cite[1.9]{FQ}) for
$S_1, S'_1, \overline S_3, \overline S'_3$. Due to the choice of the
disks bounded by the link components in $B^4$, $\overline B^4$, each
of their intersections involves one of the spheres $S_1, S'_1,
\overline S_3, \overline S'_3$. Using the embedded duals, one
resolves all extra intersection points among the spheres, getting
six smoothly embedded hyperbolic pairs. More specifically (see
\figref{double figure}) the intersections of $S_1, S'_1$ with $S_2, S'_2$ are
resolved by adding parallel copies of $\overline S_1,\overline S'_1$
to $S_2, S'_2$. Similarly, the intersections of $\overline S_2,
\overline S'_2$ with $\overline S_3, \overline S'_3$ are eliminated
by adding parallel copies of $S_3, S'_3$ to $\overline S_2,
\overline S'_2$. Surgering out the resulting embedded hyperbolic
pairs, one gets a smooth $4$-manifold $M$ homotopy equivalent to the
double cover $\Gamma$ of ${\vee}^3 S^1$.

The double cover in the general case (when the link $\Bor$ is further
Bing doubled and ramified) is constructed analogously. After
additional Bing doubling, the link can still be changed into the
unlink by intersecting a single pair of components --- there is
actually more freedom in choosing the pair of components to
intersect. The proof goes through without significant changes also
for ramified Bing doubles (when one takes parallel copies of the
components before Bing doubling them). Any such link is obtained
from the Hopf link $H=(a,b)$ by an iterated application of taking
parallel copies and Bing doubling, so that the resulting link has
trivial linking numbers. During the first iteration, at least one of
the components of $H$, say $a$, and all of its parallel copies, are
going to be Bing doubled. Label the other component, $b$,  by $1$,
and the Bing doubles of $a$ and of its parallel copies by $2$ and
$3$. (That is, the two components of each such Bing double are
labeled by $2$ and $3$ respectively --- pick any of the two possible
labelings. Note that, in particular the resulting link becomes the
unlink if one removes all components with the label $i$, for any
given $i\in \{1,2,3\}$.) During the following iterations, the
parallel copies and Bing doubles inherit the label of the component
to which the operation is applied. In this general case, construct
the double cover as above. In one lift of the $4-$ball, one only has
intersections of the components labeled by $1$ and $2$. In the other
lift, there are only intersections involving the labels $2$, $3$.
The rest of the proof is identical to the above, with individual
link components and spheres labeled by $i$ replaced with the
collections of all link components and spheres labeled by $i$,
$i=1,2,3$. \endproof

\begin{figure}[t]
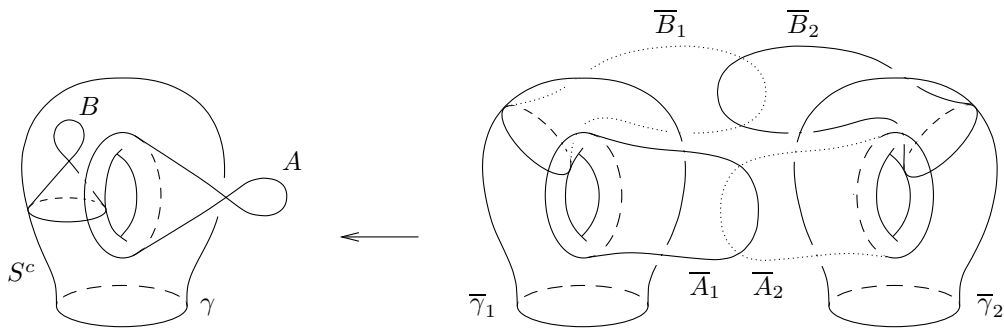
\anchor{figure double cover}
\vspace{.1cm} \small \cover \caption{The double cover, corresponding
to the homomorhism ${\pi}_1 S^c\rightarrow {\mathbb{Z}}/2$ which
sends both double point loops to the non-trivial element. The lifts
$\overline{\gamma}_1$, $\overline{\gamma}_2$ bound (non-equivariant)
disjoint embedded disks --- surgeries along the caps $\overline{A}_1$
and $\overline{B}_2$ respectively. (The caps which are not used are
drawn dotted.)} \label{figure double cover}
\end{figure}

A different class of canonical problems --- for the disk embedding
conjecture --- is provided by capped gropes, see \cite[2.1]{FQ}.
Recall that the disk embedding conjecture is the lemma underlying
the proofs of the surgery and of the s-cobordism theorems for good
groups (in fact it is equivalent to both of these theorems for any
fundamental group). Capped gropes are thickenings of certain special
$2$-complexes. They provide a class of canonical problems in the
sense that they may be found in the setting of the disk embedding
conjecture \cite[5.1]{FQ} and therefore finding a flat embedded disk
in them, with a given boundary, is equivalent to the embedding
problem in the general case. Here we present the analogue of Theorem
\ref{double cover} in this context.

Consider the simplest case (which however captures the point of the
argument): a capped surface of genus one and with just one
self-intersection point for each cap. Recall the definition of a
capped surface: start with a surface $S$ with one boundary component
$\gamma$. A capped surface $S^c$ is obtained by attaching disks to a
symplectic basis of curves in $S$. Finally, intersections are
introduced among the caps (only self-intersections of the caps, in
the current example). The interiors of the caps are disjoint from
the base surface $S$. Abusing the notation, we use $S^c$ also to
denote a special ``untwisted'' $4$-dimensional thickening of this
$2$-complex, see \cite{FQ}. The capped surface $S^c$ is homotopy
equivalent to the wedge of two circles, in particular ${\pi}_1(S^c)$
is the free group on two generators (a ``preferred'' set of
generators is given by the double point loops). Consider the
homomorphism ${\pi}_1(S^c)\longrightarrow {\mathbb Z}/2$, sending
each preferred generator to the non-trivial element, and consider
the corresponding double cover $\widetilde{S}^c$. It is an exercise
in elementary topology to check that $\widetilde{S}^c$ is given by
two copies of the capped surface of genus one, say with caps
$\overline A_1, \overline B_1$ and $\overline A_2, \overline B_2$
respectively, such that $\overline A_1$ intersects $\overline A_2$
in two points, and similarly $\overline B_1$ intersects $\overline
B_2$, \figref{figure double cover}. Note that no cap in the
cover has self-intersections, and \figref{figure double cover}
shows that the two lifts $\overline{\gamma}_1$,
$\overline{\gamma}_2$ of the boundary curve $\gamma$ bound disjoint
embedded disks in (the thickening of) $\widetilde S^c$.

Consider the general genus one case: there may be self-intersections
of each cap, and also intersections between the two caps. A slight
variation is necessary in this case: consider the homomorphism
${\pi}_1 S^c\rightarrow {\mathbb{Z}}/2$ which sends all double point
loops, corresponding to self-intersections, to $1$, and the double
point loops corresponding to the intersections of dual caps, to $0$.
The same choice of caps as in the case above works here: the caps
$\overline A_1$, $\overline B_2$ are still embedded and disjoint
(compare with the notation in \figref{figure double cover}),
since the intersections of the dual caps lift to intersections
between $\overline A_1$ and $\overline B_1$, and also to
intersections between $\overline A_2$, $\overline B_2$.

Now consider the general case, given a capped grope $(G^c,{\gamma})$
of height $n\geq 1$, whose bottom stage surface $S$ has genus $g$.
$G^c$ is obtained from $S$ by attaching capped gropes of height
$n-1$ along a symplectic basis of curves  $\{ {\alpha}_i, {\beta}_i
\}$, $i=1\,\ldots, g$ in $S$. Divide these capped gropes of height
$n-1$ into two collections, $A$ and $B$, according to whether they
are attached to one of the curves ${\alpha}_i$, or one of the
${\beta}_i$, respectively. In particular, all caps of $G^c$ are
labeled $A$ or $B$. ${\pi}_1 G^c$ is a free groups generated by the
double point loops, and generalizing the construction above,
consider the homomorphism ${\pi}_1 G^c\rightarrow {\mathbb{Z}}/2$
which sends all double point loops, corresponding to $A-A$ or $B-B$
intersections, to $1$, and the double point loops corresponding to
$A-B$ intersections, to $0$. The double cover $\widetilde G^c$
consists of two capped gropes, and neither of them has any $A-A$ or
$B-B$ intersections, since none of the double point loops of this
type in $G^c$ lift to the cover. (Here the two lifts of each cap of
$G^c$ in the cover inherit the label, $A$ or $B$, of the cap.)
Moreover, no $A-$caps of one of the gropes in the cover intersect
$B-$caps of the other one, since each double point loop of type
$A-B$ in $G^c$ lifts to a closed loop in the cover. Now the two
lifts ${\gamma}_1, {\gamma}_2$ of ${\gamma}$ bound disjoint embedded
disks: surgery along the $A-$caps in one grope, and surgery along
the $B-$caps in the other one. Thus we have proved:

\begin{lemma}
Let $(G^c, {\gamma})$ be a capped grope of height $\geq 1$, with the
attaching curve $\gamma$. Then there exists a double cover
$\overline{G}^c\longrightarrow G^c$ such that both lifts
$\overline{\gamma}_1$, $\overline{\gamma}_2$ bound disjoint smooth
disks in $\overline{G}^c$.
\end{lemma}

Instead of using the lifts of different caps to surger the gropes in
the cover, one could also use all caps, together with the operation
of contraction/pushoff (\cite{FQ}, Chapter 2.3). For example, in the
cover shown in \figref{figure double cover} the capped surfaces
can be contracted, and then the $\overline{A}_1$ --- $\overline{A}_2$
intersections are pushed off one contracted surface, while the
$\overline{B}_1$ --- $\overline{B}_2$ intersections are pushed off
the other one, to get disjoint embedded disks.

We note that the idea used here is different from the usual strategy
for the proof of the disk embedding theorem. Rather than trying to
improve the intersections between all caps, we pick certain ``good''
caps, sufficient for surgering the surface into a disk, and discard
the rest of the caps.

\vspace{.3cm}

\section{Involutions and fundamental domains.}
\label{undoubling}

\vspace{.3cm}

In this section we start the analysis of the existence of free
involutions on the family of $4$-manifolds constructed in Theorem
\ref{double cover}. This provides an approach to solving the
canonical surgery problems. Conversely, an obstruction to the
existence of such involutions would be an obstruction to surgery or
to the s-cobordism theorem for free groups. We will point out, in
particular, that the analogous problem in the closed case has a
simple solution.

The argument used in our analysis is familiar to the experts in the
$A,B$-slice problem. Stronger results are available in a related
context (cf \cite{FL}, \cite{K}) although a translation to this
setting is not immediate. We include this discussion to show a
connection of our new reformulation of the surgery conjecture with
the previous developments in the field, and to illustrate the flavor
of the problem.

Suppose there exists a required free involution on $N$ (we use the
notations of theorem \ref{double cover}). Then the quotient is a
$4$-manifold $M$ homotopy equivalent to $\vee^3 S^1$ and with
$\partial M\cong {\mathcal{S}}^0(\Wh(\Bor))$. Consider a homotopy
equivalence $M\longrightarrow \vee^3 S^1$ and lift it to a homotopy
equivalence $f$ of double covers, $f\co N\longrightarrow {\Gamma}$.
Choose $p_1,p_2,p_3 \in \vee^3 S^1$, one point in the interior of
each circle, and let $X_i=f^{-1}(p_i)$ (the transversality is
provided by \cite{Quinn}). Consider the preimages $q_i,
q'_i\in{\Gamma}$ of $p_i$, $i=1,2,3$, Figures \fref{graphcover},
\fref{equivalence}. Let $Y_i=f^{-1}(q_i)$, $Y'_i=f^{-1}(q'_i)$.
Denote $X=\cup X_i$, $Y=\cup_i Y_i$, $Y'=\cup Y'_i$. Note that due
to the ${\mathbb Z}/2$-equivariance of $f$, $Y_i$ is diffeomorphic
to $Y'_i$, for each $i$. Using the standard surgery arguments, one
may assume that each $X_i$ (and $Y_i$, $Y'_i)$ is connected. Denote
the two connected components of $N\smallsetminus(Y\cup Y')$ by $N_0$
and $N_1$.

We will consider homotopy equivalences $f$ such that $\partial
Y_i=Y_i\cap \partial N\cong\partial Y'_i= Y'_i\cap
\partial N$ is a torus, for each $i$.
We need a precise description of these tori in $\partial N$. The
manifold ${\mathcal{S}}^0(\Wh(\Bor))$ has the following convenient
description. Let $\Bor'$ denote an untwisted parallel copy of the
Borromean rings $\Bor$. Then ${\mathcal{S}}^0(\Wh(\Bor))$ is obtained
from $S^3$ by cutting out tubular neighborhoods of $\Bor\cup \Bor'$
and identifying the corresponding boundary tori $T_i$, $T'_i$,
exchanging the meridian and the longitude, $i=1,2,3$. (The proof
follows from lemma \ref{calc lemma}.) Abusing the notation, we
denote the resulting tori in ${\mathcal{S}}^0(\Wh(\Bor))$ by $T_i$
again. Start with a map ${\mathcal{S}}^0(\Wh(\Bor))\longrightarrow
\vee^3 S^1$ with these given point inverses $T_i$, and consider its
lift $f_{\partial}\co \partial N\longrightarrow {\Gamma}$. We will
consider homotopy equivalences $f$ which are extensions of
$f_{\partial}$.

It follows from the above, in particular, that $\partial N$ is
obtained from two copies of $S^3\smallsetminus (\Bor\cup \Bor')$ by
identifying the boundary torus of a regular neighborhood of a
component $l_i$ of $\Bor$ in one copy of $S^3$ with the boundary
torus for the component $l'_i$ in the other copy of $S^3$, via the
diffeomorphism exchanging the meridian and the longitude, for each
$i$.

\begin{figure}[t]\anchor{equivalence}
\vspace{.4cm} \small \equivalence \nocolon\caption{} \label{equivalence}
\end{figure}

\begin{defi}
In the general context (not assuming the existence of an involution
on $N$) we say that a homotopy equivalence $f\co
N\longrightarrow{\Gamma}$ is {\em weakly equivariant} if the
restriction of $f$ to $\partial N$ is equivariant with respect to
the obvious ${\mathbb Z}/2$ action,  and there is a diffeomorphism
from $Y_i=f^{-1}(q_i)$ to $Y'_i=f^{-1}(q'_i)$ for $i=1,2,3$,
extending the diffeomorhism of their boundary tori $\partial
Y_i\longrightarrow \partial Y'_i$, given by the involution on
$\partial N$.
\end{defi}

The existence of a weakly equivariant homotopy equivalence is a
necessary condition for the existence of an involution on $N$.
(Given an involution, both $Y$, $Y'$ are diffeomorphic to the point
inverses $X$ in the quotient.) Note that this definition does not
require that $N_0$ is homeomorphic to $N_1$, and it does not impose
any equivariance conditions on the inclusions of $Y_i, Y'_i$ into
$N_0, N_1$. However one has that $\partial N_0$, $\partial N_1$ are
diffeomorphic, and are obtained from $S^3$ by cutting out a
neighborhood of  $\Bor\cup \Bor'$ and gluing in $Y, Y'$ where the
attaching maps of $Y_i, Y'_i$ differ by the diffeomorphism of the
torus exchanging the meridian and the longitude.

Note that in the analogous context: given a homotopy equivalence $f$
from a {\em closed} $4$-manifold $N$ to a graph, the point inverses
may be arranged to be $3$-spheres \cite{KL}, up to an s-cobordism of
$N$. (Changing $N$ by an s-cobordism is fine for the applications to
the surgery conjecture.) Similarly, if $N$ has boundary but
$\partial Y_i, \partial Y'_i$ are {\em $2-$spheres} then there
exists a homotopy equivalence $f$ such that the point inverses are
$3-$balls. In the problem we consider here $\partial Y_i, \partial
Y'_i$ are {\em tori}, and we will now show that a similar, naive,
guess that there is a homotopy equivalence with $Y_i\cong Y'_i=$
solid torus is not realized. (The complexity of this problem is
precisely in the interrelation between the homotopy type and the
boundary of the manifolds $M$, $N$.)

\begin{lemma} \label{torus} 
There does not exist a weakly equivariant homotopy equivalence $f\co
N\longrightarrow {\Gamma}$ such that $Y_i=f^{-1}(q_i)$ is an integer
homology $S^1\times D^2$ for each $i$.
\end{lemma}

Note the similarity with the {\em A-B slice problem} introduced in
\cite{F1}, see also \cite{FL}, \cite{K}. Assuming the existence of
$M$ as above, it is shown in \cite{F1} that the compactification of
the universal cover $\widetilde M$ is the $4$-ball. The group of
covering transformations (the free group on three generators) acts
on $D^4$ with a prescribed action on the boundary. Roughly speaking,
this approach to finding an obstruction to surgery is in eliminating
the possibilities for the fundamental domains for such actions. Our
present approach is in terms of the fundamental domains in the
double cover of $M$, and in terms of the closely related analysis of
the point inverses. Lemma \ref{torus} eliminates just the most basic
possibility for the point inverses $Y, Y'$. Note that this is not
solely a program for finding an {\em obstruction} to surgery:
conversely, solving the problem in the affirmative would construct
an involution on $N$.

\proof First assume that $Y_i, Y'_i$ are
diffeomorphic to $S^1\times D^2$ for each $i$, and consider the
Mayer-Vietoris sequence for the decomposition $N=N_0\cup_{Y,Y'}N_1$.
The sequence splits, and we have $$ 0\longrightarrow H_2(Y\cup
Y')\longrightarrow H_2(N_0)\oplus H_2(N_1)\longrightarrow 0,$$ $$
0\longrightarrow H_1(Y\cup Y')\longrightarrow H_1(N_0)\oplus
H_1(N_1)\longrightarrow 0
$$

(all homology groups are considered with the integer coefficients).
It follows that $H_2(N_0)=H_2(N_1)=0$. Since $H_1(Y\cup Y')\cong
{\mathbb Z}^6$, the rank of one of the groups $H_1(N_0)$, $H_1(N_1)$
is $\geq 3$ --- suppose this condition holds for $N_0$. It also
follows from the sequence above that the homomorphism $H_1(\partial
N_0)\longrightarrow H_1(N_0)$ is onto.

On the other hand, $\partial N_0$ is obtained from
$S^3\smallsetminus(\Bor\cup \Bor')$ by gluing in $Y\cup Y'$, where the
attaching map of $\partial Y_i$ differs from the attaching map of
$\partial Y'_i$ by the diffeomorphism of the torus exchanging the
meridian and the longitude. Let $(p_i,q_i)$ be the slope of the
curve in $T_i$ which bounds in the solid torus $Y_i$. The
coordinates are given by (meridian, longitude) in $T_i$ which is
considered as the boundary of a tubular neighborhood of a component
of $\Bor$ in $S^3$. Then $(q_i,p_i)$ is the slope in $T'_i$ which
bounds in $Y'_i$. That is, $\partial N_0$ is the Dehn surgery on
$S^3$ along $\Bor\cup \Bor'$, with the surgery coefficients
$(p_i,q_i)$ for the component $l_i$ of $\Bor$ and $(q_i, p_i)$ for
the component $l'_i$ of $\Bor'$, $i=1,2,3$.

It follows from this description that, considering various
possibilities for the pairs $(p_i,q_i)$, the maximal possible rank
of $H_1(\partial N_0)$ is $3$. Combining this with the facts that
$rk(H_1(N_0))\geq 3$ and the map $H_1(\partial N_0)\longrightarrow
H_1(N_0)$ is onto, we conclude that $rk(H_1(\partial
N_0))=rk(H_1(N_0))=3$. Moreover, if any of the pairs $(p_i,q_i)$ is
not equal to $(1,0)$ or $(0,1)$ then the rank of $H_1(\partial N_0)$
is less than $3$. Therefore assume that each pair $(p_i,q_i)$ is
equal to either $(1,0)$ or $(0,1)$. In each of these cases $\partial
N_0={\mathcal S}^0 \Bor$, the zero-framed surgery on the Borromean
rings, and ${\pi}_1 (\partial N_0)$ is abelian.

Since $H_1(\partial N_0)\longrightarrow H_1(N_0)$ is an isomorphism
and $H_2(\partial N_0)\longrightarrow H_2(N_0)$ is onto, by
Stallings theorem \cite{S} the inclusion $\partial
N_0\hookrightarrow N_0$ induces an isomorphism on
${\pi}_1/{\pi}_1^k$ for all $k$. Here ${\pi}^k$ denotes the $k$th
term of the lower central series of a group $\pi$. Another
application of Stallings theorem, to an inclusion $\vee^3
S^1\hookrightarrow N_0$, implies that ${\pi}_1(N_0)/{\pi}_1(N_0)^k$
is isomorphic to $F_3/F_3^k$ for all $k$, where $F_3$ is the free
group on three generators. This is a contradiction since
${\pi}_1(\partial N_0)$ is abelian.

In the the general case, $Y_i$ is an integer homology $S^1\times
D^2$, for each $i$. Consider the degree one maps $Y_i\longrightarrow
S^1\times D^2, Y'_i\longrightarrow S^1\times D^2$, which are
diffeomorphisms of the boundaries. Gluing these maps with the
identity on $S^3\smallsetminus (\Bor\cup \Bor')$, one has a map
$\partial N_0\longrightarrow {\mathcal S}^0 \Bor$ inducing an
isomorphism on homology. By Stallings theorem, this map induces an
isomorphism of nilpotent quotients of the fundamental groups, in
particular ${\pi}_1(\partial N_0)$ is abelian. The Mayer-Vietoris
calculation for the decomposition $N\smallsetminus (Y\cup
Y')=N_0\cup N_1$ is still valid, showing that
${\pi}_1(N_0)/{\pi}_1(N_0)^k\cong F_3/F_3^k$ for all $k$ and giving
a contradiction as above. \qed

\Addressesr

\end{document}

%% file: agtout.tex
\def\ifplaintex{\expandafter\ifx\csname documentclass\endcsname\relax}

\def\gtp{{\mathsurround=0pt\it $\cal G\mskip-2mu$eometry \&\ 
$\cal T\!\!$opology $\cal P\!$ublications}}  % GT publications

\def\Addressesr{\bigskip
{\small \parskip 0pt \leftskip 0pt \rightskip 0pt plus 1fil \def\\{\par}
\sl\theaddress\par
\medskip
\rm Email:\stdspace\tt\theemail\hfill\rm Received:\qua\receiveddate \par}}

\def\recd{{\small Received:\qua\receiveddate\ifx\reviseddate\relax
\else\qquad Revised:\qua\reviseddate\fi\par}} 

\def\lognumber#1{\def\thelognumber{#1}}
\def\volumenumber#1{\def\thevolumenumber{#1}}
\def\volumeyear#1{\def\thevolumeyear{#1}}
\def\papernumber#1{\def\thepapernumber{#1}}
\def\pagenumbers#1#2{\def\startpage{#1}\def\finishpage{#2}}
\def\published#1{\def\publishdate{#1}}

\def\received#1{\def\receiveddate{#1}}

\def\accepted#1{\def\accepteddate{#1}}

\long\def\asciiabstract#1{\long\def\theasciiabstract{#1}}

\let\\\par\let\thelognumber\relax\let\thevolumenumber\relax
\let\thepapernumber\relax\let\thevolumeyear\relax\let\startpage\relax
\let\finishpage\relax\let\publishdate\relax\let\receiveddate\relax
\let\reviseddate\relax\let\accepteddate\relax\let\theasciititle\relax
\let\theasciiauthors\relax
\let\theasciiabstract\relax

\let\theasciiemail\relax

\ifplaintex
\font\logobig=cmssbx10 scaled 3836
\font\logomed=cmssbx10 scaled 2557
\else
\font\logobig=cmssbx10 scaled 4200
\font\logomed=cmssbx10 scaled 2800
\fi

\long\def\makeagttitle{   %%% start of definition of \makeagttitle
\count0=\startpage
\agt\hfill      %   Journal title (top left) 
\hbox to 45truept{\vbox to 0pt{\vglue -13truept{\logomed A\kern -.37em{\logobig 
T}\kern -.38em G}\vss}\hss}
\break
{\small Volume \thevolumenumber\ (\thevolumeyear)
\startpage--\finishpage\nl
Published: \publishdate}

\vglue .25truein

{\parskip=0pt\leftskip 0pt plus
1fil\def\\{\par\smallskip}{\Large\bf\thetitle}\par\medskip} \vglue
0.05truein

{\parskip=0pt\leftskip 0pt plus 1fil\def\\{\par}{\sc\theauthors}
\par\medskip}%
 
\vglue 0.03truein 

{\small\leftskip 25truept\rightskip 25truept{\bf Abstract}\stdspace\theabstract

{\bf AMS Classification}\stdspace\theprimaryclass
\ifx\thesecondaryclass\relax\else; \thesecondaryclass\fi\par
{\bf Keywords}\stdspace \thekeywords\par}\vglue 7truept

}   %%%% end of definition of \makeagttitle

\ifplaintex
\font\phead=cmsl9 scaled 950
\font\pnum=cmbx10 scaled 913
\font\pfoot=cmsl9 scaled 950
\headline{\vbox to 0pt{\vskip -4.5mm\line{\small\phead\ifnum
\count0=\startpage ISSN 1472-2739 (on-line) 1472-2747 (printed)
\hfill {\pnum\folio}\else\ifodd\count0\def\\{ }% 
\ifx\theshorttitle\relax\thetitle\else\theshorttitle\fi\hfill{\pnum\folio}
\else\def\\{ and }{\pnum\folio}\hfill\ifx\theshortauthors\relax\theauthors
\else\theshortauthors\fi\fi\fi}\vss}}
\footline{\vbox to 0pt{\vglue 0mm\line{\small\pfoot\ifnum\count0=\startpage
\copyright\ \gtp\hfill\else
\agt, Volume \thevolumenumber\ (\thevolumeyear)\hfill\fi}\vss}}
\else
\font=cmsl9 scaled 1050
\font\lnum=cmbx10 
\font=cmsl9 scaled 1050
\makeatletter
\def\@oddhead{{\small\lhead\ifnum\count0=\startpage ISSN 1472-2739 
(on-line) 1472-2747 (printed)\hfill {\lnum\number\count0}\else\ifodd\count0
\def\\{ }\ifx\theshorttitle\relax \thetitle \else\theshorttitle\fi\hfill
{\lnum\number\count0}\else\def\\{ and }{\lnum\number\count0}
\hfill\ifx\theshortauthors\relax 
\theauthors\else\theshortauthors\fi\fi\fi}}\def\@evenhead{\@oddhead}
\def\@oddfoot{\small\lfoot\ifnum\count0=\startpage\copyright\ \gtp\hfill\else
\agt, Volume \thevolumenumber\ (\thevolumeyear)\hfill\fi}
\def\@evenfoot{\@oddfoot}
\makeatother
\fi
\let\maketitlepage\makeagttitle
\let\makeshorttitle\maketitlepage
\let\maketitle\maketitlepage

\newwrite\gtoutfile
\long\gdef\makeheadfile{  %%% start of definition of \makeheadfile
{\def\\{, }\def\s{ }
\immediate\openout\gtoutfile head.xxx
\immediate\write\gtoutfile{Proxy-for: \ifx\theasciiauthors\relax
\theauthors\else\theasciiauthors\fi\s<\ifx\theasciiemail\relax\theemail\else\theasciiemail\fi>}
\immediate\write\gtoutfile{\noexpand\\}
\immediate\write\gtoutfile{Authors: \ifx\theasciiauthors\relax
\theauthors\else\theasciiauthors\fi}
{\def\\{ }\immediate\write\gtoutfile{Title: \ifx\theasciititle\relax
\thetitle\else\theasciititle\fi}}
\immediate\write\gtoutfile{Subj-class: GT or SG, GR etc}
\immediate\write\gtoutfile{MSC-class: \theprimaryclass\ifx\thesecondaryclass\relax\else, \thesecondaryclass\fi}
\immediate\write\gtoutfile{Journal-ref: Algebr. Geom. Topol. \thevolumenumber\s
(\thevolumeyear) \startpage-\finishpage}
\immediate\write\gtoutfile{Comments: Published by Algebraic and
Geometric Topology at}
\immediate\write\gtoutfile{\s\s\s  http://www.maths.warwick.ac.uk/agt/AGTVol\thevolumenumber/agt-\thevolumenumber-\thepapernumber.abs.html}
\immediate\write\gtoutfile{\noexpand\\}
\immediate\write\gtoutfile{}
\ifx\theasciiabstract\relax
\immediate\write\gtoutfile{\theabstract}\else
\immediate\write\gtoutfile{\theasciiabstract}\fi
\immediate\write\gtoutfile{}
\immediate\write\gtoutfile{\noexpand\\}
\immediate\write\gtoutfile{}
\immediate\closeout\gtoutfile}}  %%% end of definition of \makeheadfile

\def\maketitlepage{\makeagttitle\makeheadfile}
\let\makeshorttitle\maketitlepage
\let\maketitle\maketitlepage